\documentclass[11pt]{article}
 \usepackage{amsmath, amsthm, amsfonts, amssymb}

\hoffset= - 1.1 cm \voffset=  0.2 cm \numberwithin{equation}{section}

\newtheorem{theorem}{Theorem}[section]
\newtheorem{lemma}[theorem]{Lemma}
\newtheorem{proposition}[theorem]{Proposition}
\newtheorem{corollary}[theorem]{Corollary}
\newtheorem{remark}[theorem]{Remark}
\newtheorem{definition}[theorem]{Definition}
\newtheorem{example}[theorem]{Examples}
\newtheorem{hypothesis}[theorem]{Hypothesis}

\def\qed{\hfill \hbox{\hskip 6pt\vrule
width6pt height6pt depth1pt  \hskip1pt}
\smallskip}

\newcommand{\bth}{\begin{theorem}}
\renewcommand{\eth}{\end{theorem}}
\newcommand{\bpr}{\begin{proposition}}
\newcommand{\epr}{\end{proposition}}
\newcommand{\bco}{\begin{corollary}}
\newcommand{\eco}{\end{corollary}}
\newcommand{\ble}{\begin{lemma}}
\newcommand{\ele}{\end{lemma}}
\newcommand{\bpf}{\begin{proof}}
\newcommand{\epf}{\end{proof}}
\newcommand{\bex}{\begin{example}}
\newcommand{\eex}{\end{example}}
\newcommand{\bdf}{\begin{definition}}
\newcommand{\edf}{\end{definition}}
\newcommand{\bre}{\begin{remark}}
\newcommand{\ere}{\end{remark}}
\newcommand{\bhy}{\begin{hypothesis}}
\newcommand{\ehy}{\end{hypothesis}}

\newcommand{\bal} {\begin{aligned}}
\newcommand{\eal}{\end{aligned}}
\newcommand{\ben}{\begin{enumerate}}
\newcommand{\een}{\end{enumerate}}
\newcommand{\bit}{\begin{itemize}}
\newcommand{\eit}{\end{itemize}}

\def\P{{\mathbb P}}
\def\R{{\mathbb R}}
\def\C{{\mathbb C}}
\def\E{{\mathbb E }}

\def\N{{\mathbb N}}

\def\lan{{\langle}}
\def\ran{{\rangle}}
\def\hh{{\vskip 1mm} \noindent}

\begin{document}

\begin{center} { \Huge Densities for  Ornstein-Uhlenbeck processes with
jumps}
\end{center}

%\date{dd}

%\maketitle
\ \

\begin{center}
11 april 2008
\end{center}

% {\vskip 1.0 cm}

\begin{center}

 Enrico Priola \footnote
 { \noindent \! \! \! Supported by  the Italian
 National Project
 MURST ``Equazioni di Kolmogorov'' and
 by the Polish Ministry of Science and Education project 1PO 3A 034
29 ``Stochastic evolution equations with L\'evy noise''.}

\vspace{ 2 mm } {\small  \it Dipartimento di Matematica,
Universit\`a di Torino, \par
  via Carlo Alberto 10,  \ 10123, \ Torino, Italy. \par
 e-mail \ \ enrico.priola@unito.it }
\\ \ \\
Jerzy Zabczyk   \footnote { \noindent \! \!  Supported by the Polish
Ministry of Science and Education project 1PO 3A 034 29 ``Stochastic
evolution equations with \! \! L\'evy noise''.}

\vspace{ 2 mm} {\small  \it Instytut Matematyczny,   Polskiej
Akademii Nauk,
\par ul. Sniadeckich 8, 00-950, \ Warszawa, Poland.\par
 e-mail \ \ zabczyk@impan.gov.pl }
\end{center}

{\vskip 7mm }
 \noindent {\bf Mathematics  Subject Classification (2000):} \
 60H10,
% SDES
60J75,
% jump processes
 47D07.
% markov semigroups
 \ \par \ \par

\noindent {\bf Key words:}  Ornstein-Uhlenbeck processes, absolute
continuity, L\'evy processes.

\vspace{2.5 mm}

\noindent {\bf Abstract:}    We consider an Ornstein-Uhlenbeck
process with values in $\R^n$ driven by a L\'evy process $(Z_t)$
taking values in $\R^d$ with $d $ {\it possibly smaller} than $n$.
The L\'evy noise can have a degenerate or even vanishing  Gaussian
component.
 Under a
 controllability rank condition and a mild assumption on the L\'evy
 measure of $(Z_t)$, we prove that
 the law of the Ornstein-Uhlenbeck process
 at any time $t>0$ has a density on $\R^n$.
   Moreover,  when the L\'evy process  is
 of $\alpha$-stable type, $\alpha \in (0,2)$,
 %under the controllability condition,
     we show that  such
 density is a $C^{\infty}$-function.

{\vskip 1cm}
\section {\textbf{ Introduction and statement of the  main results}}

  We study  absolute continuity of the laws
   of a $n$-dimensional Ornstein-Uhlenbeck
   process $(X_t^x),$ which   solves the stochastic differential
   equation
\begin{equation}\label{1}
 dX_t = A X_t dt + B dZ_t,\;\; X_0 =x \in \R^n.
\end{equation}
Here $(Z_t)$ is a given L\'evy process, with values in $\R^d $,
defined on
 some stochastic basis $(\Omega,  {\cal F}, \! ({\cal F}_t)_{t \ge
 0},
  \P )$.
The dimension   $d$ might be different and also {\it smaller} than
$n$. Let us recall that  $(Z_t)$
 is a stochastic process having
   independent, time-homogeneous  increments and
    {\emph  c\`adl\`ag } trajectories,
     starting from $0$ (see \cite{A}). Moreover $A$ is a
      real $n \times n$
   matrix and $B$ a real $n \times d$  matrix.

Ornstein-Uhlenbeck processes appear in many areas of science, for
instance  in physics (see \cite{Garbaczewski} and the references
   therein) and in mathematical finance (see
   \cite{Barndorff-Nielsen}, \cite{Br},
     \cite{CT} and
the references
   therein). Ornstein-Uhlenbeck processes with jumps have
   recently received much  attention
     (see \cite{Za2}, \cite{SWY1},
     \cite{sato} and \cite{PZ}).

% \vspace{1 mm}

In the paper we present two main results: one   on existence of
densities of $(X_t^x),$ and the other on   the regularity of such
densities. Both theorems  assume the following (controllability)
 {\it rank condition}
 \begin{equation} \label{rank}
\! \! \! \! \! \!  \! \! \! \! \! \! \! \! \!  \;\;\; \mbox{Rank}
\, [B, AB, \ldots, A^{n-1}B]=n. \end{equation}
 Here $[B, AB, \ldots, $ $A^{n-1}B]$  denotes the $n \times nd$
 matrix,
 composed of matrices   $B, \ldots, A^{n-1}B $,
  which  corresponds to the linear
 mapping: $(u_0, \ldots, u_{n-1} ) $ $ \mapsto $
  $ B u_0 + \, \ldots$ $+ \, A^{n-1}B u_{n-1}$,
   from $\R^{nd}$ into $\R^n$.
 An interesting  example of   an Ornstein-Uhlenbeck process
   with degenerate noise satisfying the rank condition (with $d =1$ and
    $n=2$)
   %covered by our paper,
   is a solution of  the  equation
%   a generalization of the example given by Kolmogorov, i.e.,
\begin{equation}\label{Kol}\left\{
\bal  X_t^1   & =  Z_t,  \;\;\; X_0^1 = x_0^1,\\
      X_t^2  & = x_0^1 t + \int_0^t Z_s ds + x_0^2,\;\;\;
      \;\;\; t \ge 0,
      \; \; x =(x_0^1, x_0^2) \in \R^2.
\eal \right.
\end{equation}
 It is  a generalization of  a famous
 example due to Kolmogorov, in which $(Z_t)$ was a real
 Wiener process.
 In \cite{Ko} Kolmogorov  showed that
 the law of the random variable
 $(X_t^1, X_t^2)$ is absolutely
   continuous with respect to the Lebesgue measure, for any $t>0$,
    $x \in \R^2$,
   and, in fact, its density
   is a $C^{\infty}$-function on $\R^2$. This
    Gaussian
   example has also been   considered by H\"{o}rmander in
 \cite{Ho}.

% \vspace{1 mm}

\noindent When  the process $(Z_t)$ is a standard
  $d$-dimensional Wiener
process, it is known that $X_t^x$ has a density if and only if the
{ rank condition}  holds (see, e.g., \cite{DZ0} and \cite
 {DZ2}).  Moreover under \eqref{rank}
   the random variables  $X_t^x$, $t>0$, $x \in \R^n$,
 have
$C^{\infty}$-densities.   This regularity result can be easily
extended to the case when
  the L\'evy process
  $(Z_t) $ is given by  a non-degenerate
 $d$-dimensional  Wiener process with
  a drift plus an  independent pure jump process
  (see Section 2).
  Indeed, in such case,  $X_t^x$ has   two independent components,
   one of which is a  Gaussian Ornstein-Uhlenbeck process at time $t$
    having
     a $C^{\infty}$-density. Note that   convolution of two Borel
     probability
    measures   has a density as long as at least one of the two measures
     has a density  (see
\cite[Lemma 27.1]{sato}).
%%%%%%%%

However, the situation is less clear if the Gaussian component of
 $(Z_t)$ degenerates or vanishes.
% i.e., $(Z_t)$ is a pure L\'evy Jump
%process.
In this paper we consider such case.
  Indeed, we  formulate our  mild assumptions
 for absolute continuity only in terms of the L\'evy measure  of
  $(Z_t)$.
Our main first theorem is
 the following one.
\bth \label{ci} Assume the rank condition  \eqref{rank}. Assume
also that the L\'evy measure $\nu$ of $(Z_t)$ is infinite and
that there exists
 $r>0$ such that $\nu$
 restricted to the ball $\{ x \in \R^d \, : \, |x|\le r \}$
 has a density  with respect to the Lebesgue measure.
 Then, for any $t>0$ and $x \in \R^n$,
 the law of
 $X_t^x$
  is absolutely
 continuous.
\eth
 \noindent
 It also turns out (see Proposition \ref{e}) that
under the assumptions of the theorem,   the Ornstein-Uhlenbeck
process $(X_t^x)$ is strong Feller.  There are
% numerous
a number of papers dealing with the absolute continuity of laws of
degenerate
 diffusion processes
 with jumps (see \cite{B},
 \cite{NT}, \cite{Pi}, \cite{KT}
 % ,  and more recent:
 %also
% more recent papers:
  % \cite{OH}
 and  \cite{IK}).
They apply appropriate extensions of Malliavin calculus for jump
processes assuming also the well-known  H\"ormander
 condition on commutators (which becomes  the rank condition
 \eqref{rank}
  for Ornstein-Uhlenbeck processes).
 In \cite{B} it is assumed that the L\'evy measure of $(Z_t)$ has a
 sufficiently smooth density. In \cite{NT}, \cite{KT}, \cite{Pi}
 and \cite{IK}
    $\alpha$-stable type L\'evy processes $(Z_t)$ are considered.
The very weak sufficient conditions  for the absolute continuity
of the
 laws of degenerate
 diffusions with jumps,  formulated in  Theorem \ref{ci}, are new.
 Moreover, in the proof, we use    analytical methods as well as
control theoretic arguments.

% \vskip1mm

To formulate our second theorem, concerned with existence of regular
densities, we need a new hypothesis on the L\'evy measure $\nu $.
 \bhy \label{cio}
 There exist  $C >0$ and $ \alpha \in (0, 2)$,
   such that, for sufficiently small $r>0$,
    the following estimate
 holds:
\begin{equation} \label{ht}
 \int_{ \{ z \in \R^d \, : \,  |\lan z, h \ran |\le r  \}}
 \lan z, h \ran^2   \, \nu (dz) \ge
 C  \, r^{ 2 - \alpha} , \;\;\; h \in \R^d, \;\;
 \mbox{with}\; \; |h|=1.
 \end{equation}
 \ehy \noindent
  This condition was introduced  in \cite{Pi}.
  % \cite[Proposition ]{Pi} (see, in particular, Proposition
  % 1.1 and Remark 1).
  However, both   \cite{Pi} and
   \cite{IK}  prove  $C^{\infty}$-regularity of densities
   %   the law at time
   % $t$
    of solutions of  SDEs with jumps assuming
   a {\it strictly
   stronger}  version of Hypothesis \ref{cio} in which the
   integral with respect to  $\nu$ is taken over the smaller set
 $ \{ z \in \R^d \, : \,  | z |\le r  \}$.
 An interesting example of measure $\nu$ for which \eqref{ht} holds
 but the stronger hypothesis  is not verified
  is given in \cite[Remark 1]{Pi}.

Clearly,
  if  $(Z_t)$ is a $d$-dimensional
 $\alpha$-stable process
 which is rotation invariant
  (i.e., $\psi (h)
  = c_{\alpha} |h|^{\alpha}$, for $h \in \R^d$, $\alpha \in (0,2)$,
   where $c_{\alpha}$ is a positive constant)  then
   \eqref{ht} holds.
   %%%%%%%%%%%%%%%
 %   The following result generalizes \cite[Theorem 3]{Ku}.
  Thus  our next  theorem    generalizes the Kolmogorov
  regularity result concerning \eqref{Kol}
   to the case when $(Z_t) $ is a L\'evy process
   of $\alpha$-stable type.
   %  $\alpha \in (0,2)$.
\begin{theorem} \label{t} Assume the rank condition  \eqref{rank} and
   Hypothesis \ref{cio}.
  Then, at any time $t>0$, $x \in \R^n$,
   the Ornstein-Uhlenbeck process  $(X_t^x)$ has
  a $C^{\infty}$-density with all bounded derivatives.
  % at any time $t>0$, $x \in \R^n$.
  Moreover,
  %all the derivatives of the density are bounded
  %on $\R^n$
   for any $t>0$, $x \in \R^n$, $f: \R^n \to \R$ Borel and
   bounded,
\begin{eqnarray}
\\
% P_t f (x)
\nonumber \E [f(X_t^x)] & = & \frac{1}{(2 \pi)^n} \int_{\R^n } \!
\! f (e^{tA} x + y ) \Big(\!  \int_{\R^n} e^{- i \lan y, h\ran}
\exp { \! \Big(- \! \! \int_0 ^t \psi (B^* e^{s A^*}  h) ds \Big)}
dh \Big)dy
\\ \nonumber
& = & \frac{1}{(2 \pi)^n}  \! \int_{\R^n } \!  \!  \! f (  z )
\Big (  \!\int_{\R^n}  \!  \! e^{- i \lan z, h\ran} e^{  i \lan
e^{tA^*}  h, x\ran} \exp { \Big( \! -  \! \! \int_0 ^t \psi( B^*
e^{s A^*} h)  ds \Big)} dh \Big) dz.
%& f \in B_b(\R^n).
\end{eqnarray}
\end{theorem}

\section{Existence of  densities}

 Consider the Ornstein-Uhlenbeck process introduced  in \eqref{1}.
 It is well known that this is given by
\begin{equation}\label{2}
 X^x_t  =  e^{tA} x + \int_0^t e^{(t-s)A} B dZ_s = e^{tA} x + Y_t,
  \;\; t \ge 0, \; x \in \R^n,
\end{equation}
where the stochastic convolution $Y_t$ can be defined as a limit
in probability of Riemann  sums (see, for instance, \cite[Section
17]{sato} and \cite{SWY1}).

The law $\mu_t^x $ of $X_t^x$ has the characteristic function (or
Fourier transform) $\hat \mu_t^x$,
\begin{equation}\label{c}
\hat \mu_t^x (h) = e^{i \lan e^{tA}  x, h\ran} \hat \mu_t (h) =
e^{i \lan e^{tA^*}  h, x\ran}
 \,  \exp { \Big(- \int_0 ^t \psi ( B^* e^{s A^*}  h) ds \Big)}, \;\; h \in
 \R^n,
\end{equation}
where $\mu_t$ denotes the law of $Y_t$ and $\psi$ is the exponent
 of $(Z_t)$,
 $$
 \E [e^{i \lan u, Z_t\ran}] = e^{- t
\psi(u)},\, u \in \R^d.
$$
By $\lan \cdot, \cdot \ran$ and $|\cdot |$ we indicate the inner
product and the Euclidean norm in  $\R^k$, $k \in \N$,
respectively.
 Moreover $B^*$ denotes the adjoint (or transposed)
 matrix of $B.$

\noindent Recall the L\'evy-Khintchine representation for $\psi$,
\begin{equation} \label{ft1}
 \psi(s)= \frac{1}{2}\lan Q s,s  \ran   - i \lan a, s\ran
- \int_{\R^d} \Big(  e^{i \lan s,y \ran }  - 1 - \, { i \lan s,y
\ran} \, I_D \, (y) \Big ) \nu (dy), \;\;\; s \in \R^d,
\end{equation}
 where $I_{D}$ is the indicator function of the ball $D =
 \{ x \in \R^d \, : \, |x| \le 1 \}$,
   $Q$ is a symmetric $d \times d$ non-negative definite matrix, $a
\in
 \R^d$,
 and $\nu$
 is the  L\'evy measure of $(Z_t).$ Thus  $\nu$ is a $\sigma$-finite
measure on $\R^d$, such that
$$
\nu (\{ 0\})=0, \;\;\;\; \int_{\R^d} (1 \wedge |y|^2 ) \, \nu(dy)
 <\infty.
$$
The triplet $(Q, a, \nu)$ which gives (\ref{ft1})
   is unique. According to \eqref{ft1}, the process $(Z_t)$
    can be represented by the L\'evy-It\^o decomposition
    as
   \begin{equation}  \label{ztt}
 Z_t = at + R W_t + Z^{0}_t, \;\; \;\;\; t \ge 0,
   \end{equation}
 where $R$ is a $d \times d$ matrix such that $R R^* = Q$,
 $(W_t)$  is a standard $\R^d$-valued Wiener process and
 $(Z^{0}_t)$ is a L\'evy  jump process (see \cite{A}).
  The processes $(W_t)$ and $(Z^{0}_t)$
 are independent.

\vskip 1mm Let $(P_t)$ be the transition semigroup determined by
$(X_t^x)$,
 i.e.,
$$
 P_t f (x) = \E [f (X_t^x)],\;\; t \ge 0, \; x \in \R^n,
$$
 $f \in B_b (\R^n)$, where $B_b (\R^n)$ denotes the space of all real
  Borel and bounded functions on $\R^n$.
  The semigroup $(P_t)$ (or the process $(X_t^x)$)
   is called {\it strong Feller} if $P_t f $  is a continuous
   function,
  for any $t>0$ and for any $f \in B_b (\R^n)$.

\hh Applying a result due to Hawkes (see \cite{Ha}) we  show now
 that the strong Feller property for $(P_t)$ is equivalent  to the
existence of a density for the law of $X_t^x$, for any $t>0$, $x
\in \R^n$.
 This result holds for any Ornstein-Uhlenbeck
  process defined in \eqref{1}
  (without requiring
 the rank condition \eqref{rank}).
  For related results in  infinite dimensions, see \cite{RW}
  and \cite{PZ1}.
 \begin{proposition} \label{e}
The semigroup $(P_t)$ is strong Feller  if and only if, for each
$t>0$, $x \in \R^n$, the law $\mu_t^x $ of $X_t^x$ is absolutely
continuous with respect to the Lebesgue measure.
\end{proposition}
\begin{proof} Fix $t>0$ and let
 $\mu_t $ be the law of $Y_t$ (see \eqref{2}).  Since
 $\mu_t^x = \delta_{e^{tA}x} * \mu_t$ (where $\delta_a$ denotes
 the  Dirac
  measure concentrated in $a \in \R^n$)
   $\mu_t^x $ is absolutely continuous, for any  $x \in \R^n$,
    if and
  only if $\mu_t$ has the same property.

 We write, for any $f \in B_b(\R^n)$, $x \in \R^n,$
\begin{eqnarray*}
 && P_t f (x) = \int_{\R^n} f (e^{tA} x  + y) \mu_t (dy) =
\int_{\R^n} (f \circ e^{tA}) ( x  + e^{-tA} y) \mu_t (dy) =
 \\ && \int_{\R^n} (f \circ e^{tA}) ( x  + z) (e^{-tA} \circ \mu_t)
 (dz),
\end{eqnarray*}
where $(e^{-tA} \circ \mu_t)$ is the image of the probability
measure  $\mu_t$ under $e^{-tA}$. Applying \cite[Lemma 2.1]{Ha},
we know that the Markov operator  $T_t g (x) = \int_{\R^n} g ( x +
z) (e^{-tA} \circ \mu_t)
 (dz)$, $x \in \R^n$, maps Borel and bounded functions into
 continuous ones if and only if $(e^{-tA} \circ \mu_t)$ is
absolutely continuous with respect to the Lebesgue measure.
 Hence $P_t f$ is continuous, for any $f
\in B_b (\R^n)$, if and only if $(e^{-tA} \circ \mu_t)$ is
absolutely continuous. This gives the assertion, since $e^{tA} $
is an isomorphism.
\end{proof}

\bre {\em If the Ornstein-Uhlenbeck process $(X_t^x) $ has a
 density for any $x \in \R^n$, $t>0$, i.e.,
 $P_t f (x) = \int_{\R^n} f(e^{tA}x +  y ) g_t(y) (dy) $, then $P_t f $
  is
 {\it uniformly continuous on $\R^n$,} for any $f \in L^{\infty}(\R^n)$
 and $t>0$.
 To prove this, fix $t>0$, $f \in L^{\infty}(\R^n)$
   and consider a sequence
 $(g_t^k)$ of continuous functions having compact support which
 converges to $g_t$ in $L^1(\R^n)$.
 Define, for any $k \in \N$, $P_t^k f : \R^n \to
 \R$,
 $
 P_t^k f(x) = \int_{\R^n} f(e^{tA}x +  y ) g_t^k(y) (dy)
$, $x \in \R^n$. We have that  $(P_t^k f )$ converges to $P_t f$
uniformly on $\R^n$ and,  moreover, each function $P_t^k f$ is
uniformly continuous on $\R^n$. It follows that $P_t f $ is
 uniformly continuous as well.
 %%%%%%%%%%%%%%%
 } \ere
 %%%%%%%%%%
The proof of Theorem \ref{ci} requires two lemmas.  The first one
is of independent interest.
 \ble \label{fr} Assume the rank
 condition \eqref{rank}.
 % where $A \in L(\R^n)$
 %and $B \in L (\R^d, \R^n)$.
 Then there exists $T_0 >0$
 (depending on the dimension  $n$ and on the eigenvalues of $A$)
 such that for any integer   $m \ge n+1$,
 % Then, there exists $\hat t >0$
 %(depending on the dimension $n$ and on the eigenvalues of $A$)
 %such that,
  for any $0 \le s_1 < \ldots < s_m \le T_0$,
  the  linear transformations  $
 l_{s_1, \ldots, s_m} : \R^{dm} \to
 \R^n$,
 \begin{equation}  \label{d4}
l_{s_1, \ldots, s_m}  (y_1, \ldots, y_m) = \sum_{j=1}^m e^{s_j A}B
y_j, \; \; \;\;\; \mbox{are onto.} \end{equation}
% \mbox{ \it
% are onto}.
\ele
 \bpf The proof is divided into two parts.

\hh {\it I Part.} {\it We define  $T_0>0$.}

\hh Let $(\lambda_j)$ be the distinct complex eigenvalues of $A$,
$j =1,\ldots ,k$ (with $k \le n$). Consider the following complex
polynomial: $p(\lambda) = \prod_{j=1}^k (\lambda -
\lambda_j)^{n}$, $\lambda \in \C$, and the corresponding ordinary
linear differential operator $p (D)$ of order $n$,
 $$
 p(D) y (t) = \big (\prod_{j=1}^k (D - \lambda_j)^{n} \big ) y(t)
 = y^{(n)} (t) + a_{1} y^{(n-1)}(t) + \ldots + a_n,\;\;
\; t \in \R,
$$
 where $y \in C^{n}(\R) $,
  $a_i \in \C $,
  and $y^{(i)}$ denotes the $i$-derivative of $y$,
  $i = 1 , \ldots, n$.

  By a result due to Nehari (see \cite{N})
  we know, in particular, that there exists  $T_0 >0$ (depending on $n$ and
  on the  coefficients
  $a_1, \ldots, a_n$) such that {\it any } non-trivial solution $y (t)$
  to the equation $p(D) y =0$ has at most $n$ zeros on $[-T_0, T_0]$.
 By this theorem, we deduce that   the following quasi-polynomials
 \begin{equation} \label{dd}
y(t)= \sum_{j=1}^k \sum_{r=0}^{n-1} c_{rj} e^{\lambda_j t} t^r,
 \end{equation}
 which are  solutions for
   $p(D)y=0$ (see, for instance,
  \cite[Chapter 3]{BR}),
  have always at most $n$ zeros on $[-T_0, T_0]$  no matter what are
  the complex coefficients $c_{rj}$ (except the trivial case in which
    all $c_{rj}$
   are zero).

\hh {\it II Part.} {\it We prove the assertion.}

\hh
 Introduce the following linear and bounded operators
  (depending on $t>0$)
  \begin{eqnarray*}
 L_t : L^2([0,t]; \R^d) \to \R^n,\;\;\; L_t u =
  \int_0^t e^{sA}B u(s)ds, \;\;\; u \in L^2([0,t] ; \R^d).
\end{eqnarray*}
 The controllability condition is equivalent to the fact that
 each $L_t$ is onto,   $t>0$ (see, for instance, \cite[Chapter 1]{Za}).
  Hence, in particular,
  Im$(L_{T_0}) =\R^n$ ($T_0>0$ is defined in the first part of the proof).
    To prove the assertion it is enough to show that
 \begin{equation} \label{tr}
 \mbox{Im} (L_{T_0}) \subset \mbox{Im \!}(l_{s_1, \ldots, s_m})
 \end{equation}
 for any $0 \le s_1 < \ldots < s_m \le T_0$, and  $m \ge n+1$.
  We fix $m \ge n+1$ and take $(s_1, \ldots, s_m)$ with
$0 \le s_1 < \ldots < s_m \le T_0$.
  Let $v \in \R^n$, $v \not =  0$,
   be orthogonal to Im$(l_{s_1, \ldots, s_m})$. Assertion \eqref{tr}
    follows if we prove that
 \begin{equation} \label{lt}
 \lan v, L_{T_0} u \ran =0,\;\;\;\; \mbox{for any } \;  u \in L^2
 ([0,T_0];\R^d).
\end{equation}
 To this purpose,  note that
  the orthogonality of  $v$  to Im$(l_{s_1, \ldots, s_m})$ is
 equivalent to $B^* e^{s_j A^*} v =0$, for $j =1, \ldots , m$,
 i.e.,
 \begin{equation} \label{14}
  \lan B^* e^{s_j A^*} v , e_k \ran =0 ,\;\;\; j =1, \ldots ,
  m,\;\;
   k  =1, \ldots , d,
 \end{equation}
where $(e_k)$ is the canonical basis in $\R^d$. Note that each
mapping $s \mapsto  \lan B^* e^{s A^*} v , e_k \ran$, $k = 1,
\ldots, d$, is a quasi-polynomial  like  \eqref{dd}. Since
 $m \ge n+1$, condition \eqref{14} implies that
 each mapping $\lan B^* e^{s A^*} v , e_k \ran$ is identically zero
  on $[0,T_0]$ by the first part of the proof.

It follows that
$$
\lan L_{T_0} u, v \ran  =\int_0^{T_0} \lan u(s), B^* e^{sA^*} v
\ran ds =0,
$$
 for any $u\in L^2([0,T_0]; \R^d)$.
  This implies that $v$ is orthogonal to Im$(L_{T_0} )$ and so
  \eqref{tr} holds. The proof is complete.
\epf

%%%%%%%%%%%%%%

 \ble \label{basic} Let $L : \R^p \to  \R^q $, $p \ge q$, be an onto linear
 transformation. Let $\gamma$ be a probability measure on $\R^p$ having a
 density $h $ (with respect to the Lebesgue measure).
 Then the probability measure $L \circ \gamma$,  image of $\gamma
 $
under $L $, has a density on $\R^q$. \ele
 \bpf
   Since the result is clear when $p=q$, let us assume that $p>q$.
   We identify $L$ with a $q \times p $ matrix with respect to the
   canonical bases $(f_i)_{1 \le i \le p}$ in $\R^p$
    and  $(e_i)_{1 \le i \le q}$ in $\R^q$.
    Consider the transposed matrix
    $L^*$ and complete the system of vectors
    $L^* e_1, \ldots $, $L^* e_q$ with
    vectors $f_{i_1}, \ldots f_{i_{p-q}}$ in order to get a
    basis in $\R^p$.
  Define an invertible $p \times p$ matrix $S$ having the
     vectors $L^* e_1, \ldots $, $L^* e_q, $
     $f_{i_1}, \ldots f_{i_{p-q}}$ as rows. If $\pi : \R^p \to
     \R^q$ is the projection on the first $q$ coordinates, we
     have that $L = \pi \circ S $.
    Indeed, for any $ x \in \R^p$,
$$
\pi  (S x)
%= \pi ( \lan e_1, Lx \ran_{\R^q}, \ldots , \lan e_q, Lx
%\ran_{\R^q}, \lan f_{i_1}, x  \ran_{\R^p}, \ldots, \lan
%f_{i_{p-q}}, x \ran_{\R^p})
= ( \lan e_1, Lx \ran_{\R^q}, \ldots , \lan e_q, Lx \ran_{\R^q}) =
Lx.
 $$
 Fix any Borel set $B \subset \R^q$. Using also the Fubini
 theorem, we get
$$
 \int_{\R^q} I_{B} (x) (L \circ \gamma) (dx) =
 \int_{\R^p} I_{B} ( \pi (S z))  h(z) dz  =
  \frac{1}{|\mbox{det} (S)|} \int_{\R^p} I_{B} ( \pi (y) )
   h( S^{-1}(y)) dy
$$
$$
= \frac{1}{|\mbox{det} (S)|} \int_B dy_1 \ldots dy_q
\int_{\R^{p-q} }    h \circ S^{-1}( y_1 , \ldots , y_p) dy_{q+1}
\ldots dy_{p}.
$$
It follows that $L \circ \gamma $ has the density $$
 (y_1, \ldots , y_q) \mapsto
 \frac{1}{|\mbox{det} (S)|} \int_{\R^{p-q}
}    h \circ S^{-1}( y_1 , \ldots , y_p)
 dy_{q+1} \ldots dy_{p}.
$$  \epf
%%%%%%%%

\noindent \textbf {Proof of Theorem \ref{ci}.} We will
 use Lemma \ref{fr} and adapt the
 method of
  the  proof of
 \cite[Theorem 27.7]{sato}, based on
  \cite{T} and \cite{FV}.
 %using  Lemma \ref{fr}.
 Let $T_0>0$ be as
in Lemma \ref{fr}. Using Proposition \ref{e} and the semigroup
property of $(P_t)$, in order  to prove the assertion it is enough
to show that the law of $Y_t$ (see \eqref{2}) is absolutely
continuous for any $t \in (0,T_0)$.

%Note that
% by changing the L\'evy process $(Z_s)$,
%$s \in [0,t]$, with the L\'evy process $(Z_{t-s}) $, $s \in
% [0,t]$, we can assume that
% the law of $Y_t$ is the same of $\int_0^t e^{sA} B dZ_s$. Hence,
% in this proof, we set
% \begin{equation}
%Y_t := \int_0^t e^{sA} B dZ_s.
% \end{equation}
%Moreover it is not restrictive to assume that $Q=0$ in \eqref{ft1}
%(i.e., $(Z_t) $ has no a Gaussian component).
Recall that for an arbitrary  Borel measure $\gamma$ on $\R^n$, we
have the unique
 measure decomposition
 \begin{equation} \label{dec}
 \gamma = \gamma_{ac} + \gamma_{s}
 \end{equation}
 where $\gamma_{ac}$ has a density and $\gamma_{s}$ is singular
  with respect to the Lebesgue measure.

 Define, for $N \in \N$  sufficiently large, say $N \ge N_0$ with
  $1/N_0 <r$, the measure $\nu_N$ having density
  $I_{ \{ 1/ N \le |x| \le r\} }$ with respect  to $\nu$, i.e.
$$ \nu_N = \nu \, I_{ \{ 1/ N \le |x| \le r\} } \;\;
 \mbox{and } \;\; Z_t^N
= \sum_{0< s \le t, \,  \frac{1}{N} \le |{\triangle Z}_s| \le r }
\triangle Z _s ,\;\;  t \ge 0,
$$
%%%%%%%%%
(the  measure $\nu_N$ has
 density
  $I_{ \{ 1/ N \le |x| \le r\} }$ with respect  to the measure
   $\nu$ defined in  \eqref{ft1}) and
    ${\triangle Z}_s = Z _{s} -
Z_{s-}$  ($Z_{s-}= \lim_{h \to 0^-} Z_{s+h}$).
  The process $(Z_t^N)$ is  a compound Poisson process
 and its L\'evy measure is just  $\nu_N$. By the hypotheses,
  for any $N \ge N_0$,
 $\nu_N$   {\it   has a density.} Moreover
 since $\nu $ is infinite,  we have that
 $$ c_N = \nu_N (\R^d) \to \infty
 \;\;\; \mbox{ as} \,\,\,\,\, N \to \infty.
 $$
It is well known that $(Z^N_t)$ and $(Z_t - Z^N_t)$ are
independent L\'evy processes (see, for instance, \cite{A}
 or \cite[Chapter 1]{P}). It follows, in particular, that the random
variables
 \begin{equation} \label{d5}
Y_t^N =  \int_0^t e^{(t - s)A} B dZ_s^N \;\; \mbox{ and} \;\; Y_t
-Y_t^N = \int_0^t e^{(t -s)A} B d (Z - Z^N)_s \;\;  \mbox{are
independent,}
 \end{equation}
 for any $N\ge N_0$, $t>0$. Fix $t \in (0,T_0)$ and denote by $\mu$
 {\it the law of the random variable $Y_t$ and by $\mu_N$ the one of
 $Y_t^N$.}

 Since $\mu = \mu_N * \beta_N $ (where $\beta_N$ is the law of $Y_t
- Y_t^N$), we have by \eqref{dec}
$$
 \mu = (\mu_N)_{ac}
  * (\beta_N)_s + (\mu_N)_{s} * (\beta_N)_{s} +
  (\mu_N)_{ac} * (\beta_N)_{ac} + (\mu_N)_{s}
  * (\beta_N)_{ac}.
$$
By  \cite[Lemma 27.1]{sato}) we deduce that $(\mu_N)_{ac}
  * (\beta_N)_s   +
  (\mu_N)_{ac} * (\beta_N)_{ac} + (\mu_N)_{s}
  * (\beta_N)_{ac}$ is absolutely continuous and so
 $\mu_{s}  = \big ((\mu_N)_{s} * (\beta_N)_{s} \big )_s$
 and
 \begin{equation} \label{fi}
 \mu_{s} (\R^n) \le (\mu_N)_{s} * (\beta_N)_{s} \big
(\R^n) \le  (\mu_N)_s (\R^n), \;\;\; \mbox{for any}\; N \ge
 N_0.
 \end{equation}
Now we  compute $\mu_N$ which  coincides with the law
   of  $\int_0^t e^{sA} B dZ_s^N$.

First, note that
 the law of $Z_t^N$ is given by
 $$
e^{-c_N t} \delta_0 + e^{-c_N t}  \sum_{k \ge 1} \frac{(c_N
t)^k}{k!} (\tilde{\nu_N})^k, \;\; \mbox{where} \;\; c_N =
\nu_N(\R^d),\;\;
 \; \tilde{\nu_N } = \frac{\nu_N}{c_N},
$$
  $(\tilde{\nu_N})^k = \tilde \nu_N * \ldots * \tilde \nu_N $
 ($k$-times).
 Then consider a  sequence $(\xi_i)$
of independent random variables having  the same exponential law
of intensity
   $c_N$.
 Introduce another sequence $(U_i) $ of independent random
 variables (independent also of $(\xi_i)$) having  the same  law
   $\tilde \nu_N$.

  It is not difficult to check that
  % $Y_t^N$
%   has the same law of  $\int_0^t e^{sA} B dZ_s^N$.
% Hence,
 the probability measure $\mu_N$ coincides with the law of
   the following random variable:
   % We have the following identity in distribution sense
$$
 0 \cdot 1_{ \{ \xi_1 >t \} }  +
 \sum_{k \ge 1} 1_{ \{ \xi_1 + \ldots + \xi_k \le t <
   \xi_1 + \ldots + \xi_{k+1} \} }
  \Big ( e^{\xi_1 A} B U_1 + \ldots + e^{(\xi_1 + \ldots + \xi_k) A}
   B U_k \Big ).
$$
 Note that  the events $H_0 =\{ \xi_1 >t \} $,
   $ H_k =\{ \xi_1 + \ldots + \xi_k \le t <
   \xi_1 + \ldots + \xi_{k+1} \}$ are all
   disjoint, $k \ge 1$. Since, for any $f \in B_b(\R^n)$,
$$
f (\sum_{k \ge 0} X_k 1_{H_k})  =\sum_{k \ge 0} f (X_k) 1_{H_k},
$$
  where $X_0 =0$ and
   $X_k = e^{\xi_1 A} B U_1 + \ldots + e^{(\xi_1 + \ldots + \xi_k) A}
   B U_k $, $k \ge 1$, we get
\begin{align*}
 & \E f(Y_t^N) = e^{- c_N t } f(0) + R_N, \;\; \mbox{where} \;
\\
 & R_N  =  \E f \big ( \sum_{k \ge 1} 1_{ \{ \xi_1 + \ldots +
\xi_k \le t <
   \xi_1 + \ldots + \xi_{k+1} \} }
  \big ( e^{\xi_1 A} B U_1 + \ldots + e^{(\xi_1 + \ldots + \xi_k) A}
   B U_k \big ) \Big )
% \end{eqnarray*}
%\begin{eqnarray*}
\\
& =  \sum_{k \ge 1} \E f \Big(  1_{ \{ \xi_1 + \ldots + \xi_k \le
t <
   \xi_1 + \ldots + \xi_{k+1} \} }
  \big ( e^{\xi_1 A} B U_1 + \ldots + e^{(\xi_1 + \ldots + \xi_k) A}
   B U_k \big ) \Big )
\\
& = \sum_{k = 1}^{\infty} \int_{t_1 + \ldots + t_k \le t \le
  t_1 + \ldots + t_{k+1} } (c_N)^{k+1} e^{- c_N (t_1 + \ldots + t_{k+1} )}
   d t_1 \ldots dt_{k+1}
 \\ & \cdot \,  \int_{\R^{dk}} f(
 e^{t_1 A} B y_1 + \ldots + e^{(t_1 + \ldots + t_k) A}
   B y_k \big )  \tilde \nu_N (dy_1)\ldots \tilde \nu_N (dy_k)
\\
& = \sum_{k \ge 1} \int_{t_1 + \ldots + t_k \le t \le
  t_1 + \ldots + t_{k+1} } (c_N)^{k+1} e^{-c_N (t_1 + \ldots + t_{k+1} )}
   d t_1 \ldots dt_{k+1} \cdot \\
   & \cdot \int_{\R^{n}} f( y  )
   \mu_{t_1, \ldots, t_k} (dy), \;\;\; f \in B_b(\R^n),
\end{align*}
 where $ \mu_{t_1, \ldots, t_k} $ is the probability
   measure on $\R^n$ which is the image of the product
 measure $ \tilde \nu_N \times \ldots \times \tilde \nu_N$ ($k$-times)
  under the linear
 transformation $J_{t_1, \ldots, t_k}$
 (independent of $N$) acting from $\R^{dk}$ into $\R^n$,
 $$
J_{t_1, \ldots, t_k}   ( y_1 , \ldots , y_k ) = e^{t_1 A} B y_1 +
\ldots + e^{(t_1 + \ldots + t_k) A}
   B y_k,
$$
 where $y_i \in \R^d,$ $i =1, \ldots, k$.
  For any $k \ge n+1$, $t_1 \ge 0$, $t_i >0$, $i=2, \ldots , k$,
  we have
 $ 0 \le t_1 < \ldots < t_1 + \ldots + t_{k} \le  T_0  $ and
$$
J_{t_1, \ldots, t_k}  = l_{t_1, \ldots, t_1 + \ldots + t_k}
$$
(see \eqref{d4} and recall that
  $t
\in (0,T_0)$).  Applying Lemma \ref{fr}, we obtain that, for any
$k \ge n+1$,  $t_i >0$, $i=1, \ldots , k$,
 the linear transformation $J_{t_1, \ldots, t_k}$ is {\it onto}.
 Therefore, by Lemma \ref{basic},
the measure $\mu_{t_1, \ldots, t_k}$ has a density $g_{t_1,
\ldots, t_k} \in L^1(\R^n)$, for any $k \ge n+1$, $t_i
>0$, $i=1, \ldots , k$.
 Using this fact, we write
 \begin{eqnarray*}
&& \mu_N =  \mu_N^1 + \mu_N^2, \;\; \mbox{where} \;\;\;  \mu_N^1 =
e^{- c_N t } \delta_0 + \\ && +  \, \sum_{k =1}^{ n} \int_{t_1 +
\ldots + t_k < t <
  t_1 + \ldots + t_{k+1} } (c_N)^{k+1} e^{-c_N (t_1 + \ldots + t_{k+1} )}
   \mu_{t_1, \ldots, t_k} \,
    d t_1 \ldots dt_{k+1},
\end{eqnarray*}
% && \mbox{and} \;\;\; \mu_N^2  \; \mbox{has the following
%density on } \;\; \R^n:
and $ \mu_N^2 $  has the following density on $ \R^n:$
$$
 y \mapsto \sum_{k
> n} \int_{t_1 + \ldots + t_k < t <
  t_1 + \ldots + t_{k+1} } (c_N)^{k+1} e^{-c_N (t_1 + \ldots + t_{k+1} )}
   g_{t_1, \ldots, t_k } (y)
    d t_1 \ldots dt_{k+1}.
$$
 Therefore
 $$
\bal & (\mu_N)_s (\R^n) \le  \mu_N^1 (\R^n) \\ &  =  e^{- c_N t }
+ \sum_{k =1}^{ n} \int_{t_1 + \ldots + t_k < t <
  t_1 + \ldots + t_{k+1} } (c_N)^{k+1} e^{-c_N (t_1 + \ldots + t_{k+1} )}
    d t_1 \ldots dt_{k+1}
  \; \longrightarrow \,\, 0, \eal $$
 $\mbox{as}\; N  \to \infty, $
  since
$c_N \to \infty$ by hypothesis.   By \eqref{fi}, we immediately
 get that
 $
\mu_{s}  = 0. $
  This gives the assertion. The proof is complete.

\qed

\section{Proof of the $C^{\infty}$-result}

We pass now to the {\it proof of Theorem \ref{t}}. To obtain
$C^{\infty}$-regularity of the law at time $t$ of the
Ornstein-Uhlenbeck process \eqref{1} we will be estimating its
characteristic function.

\vskip2mm

We fix $t>0$. It is enough to show that the law $\mu_t$
 of $Y_t$ (see \eqref{2}) has a density $p_t \in L^1(\R^n)
  \cap C^{\infty} (\R^n)$ with all bounded derivatives.
  To this purpose, note that by \eqref{c} the
  characteristic function of $\mu_t$ is
$$
 \hat \mu_t (y) =
  \exp  \Big(- \int_0 ^t \psi ( B^* e^{s A^*}  y) ds  \Big),\;\;\;
 y \in \R^n.
$$
%%%%%%%%
 We claim that there exist $a_t$ and  $c_t >0$ such that,
 for any $y \in \R^n$, $|y| \ge 1$,
 \begin{equation} \label{df}
 \big |\exp  \Big(- \int_0 ^t \psi ( B^* e^{s A^*}  y) ds  \Big)
  \big|
 \le   c_t e^{- a_t |y|^{\alpha} }.
 \end{equation}
 %\textbf{
 This will imply in particular that $\hat \mu_t \in L^{1}
 (\R^n)$.  Then, by using the Fourier inversion formula (see
    \cite[Propositions 2.5]{sato})
 we  will get  the assertion.
 %%%%%%%%%%%%%%%%%%%%
%%%%%%%%%%
 \vskip 1mm
  It is not restrictive to assume that $Q=0$ and $a=0$ in \eqref{ft1}, i.e.,
that $(Z_t) $ has no  Gaussian component.
  For any $y \in \R^n$, we have
   $$
   \Big |\exp  \Big(- \int_0 ^t \psi ( B^* e^{s A^*}  y) ds  \Big)
  \Big|  =   \exp  \Big(- \int_0 ^t ds
   \int_{\R^d} \big (1 - \cos  (\lan B^* e^{s A^*}  y, z \ran) \big)   \nu
   (dz) \Big).
  $$
   First,  note that condition \eqref{ht} is equivalent
    to the fact that
\begin{equation} \label{ff}
\int_{ \{z \in \R^d \, :\,  |\lan z, k \ran |\le 1  \}}
 \lan z, k \ran^2   \, \nu (dz) \ge
 C  \,  |k|^{\alpha},
\end{equation}
 for sufficiently large $k \in \R^d$, say $|k| \ge c_0$. To see
 this,
 it is enough to  change in the condition \eqref{ht},
  the vector  $h$ to the  vector $k/r$.
  Fix $y \in \R^n$ with $|y|\ge 1$;
   using also the inequality $1 -  \cos (u) \ge c_1 |u|^2$, if $|u|
\le
 \pi$,  we find
  $$ \bal &  \int_0 ^t ds
   \int_{\R^d} \big (1 - \cos  (\lan B^* e^{s A^*}  y, z \ran) \big)   \nu
   (dz)
\\
& \ge c_1 \int_0 ^t ds
   \int_{ \{ z \in \R^d \, :\,  |\lan B^* e^{s A^*}  y, z \ran| \le 1 \}
} \, \lan B^* e^{s A^*}  y, z \ran ^2   \, \nu
   (dz)
\\
& \ge c_1 \int_0 ^t 1_{ \{ s \in [0,t] \, :\,  |B^* e^{s A^*} y|
\ge c_0 \} } \, ds
   \int_{ \{ z \in \R^d \, :\,  |\lan B^* e^{s A^*}  y, z \ran| \le 1 \}
} \, \lan B^* e^{s A^*}  y, z \ran ^2   \, \nu
   (dz)
\\
& \ge c_1 C \, \int_0 ^t  1_{ \{ s \in [0,t] \, :\,  |B^* e^{s
A^*} y| \ge c_0 \} } \,  | B^* e^{s A^*}  y|^{\alpha} ds. \eal $$
%Set $ \displaystyle{ M_t = \sup_{s \in [0,t],\, |h|\le 1, \, h \in
%\R^n}
%    | B^* e^{s A^*}  h |}$;
Set  $  M_t = \sup \{ | B^* e^{s A^*}  h | \, :\,   s \in [0,t],\,
|h|\le 1, \, h \in \R^n \}$;
  since $ \Big |
\frac{ B^* e^{s A^*}  y }
 {|y| \, M_t } \Big| \le 1$, $s \in [0,t]$, we get
$$
\bal & c_1 C \, \int_0 ^t  1_{ \{ s \in [0,t] \, :\,  |B ^* e^{s
A^*} y| \ge c_0 \} } \,  | B^* e^{s A^*}  y|^{\alpha} ds
\\
 & \ge   c_1 C \,  |y|^{\alpha} \, M_t ^{\alpha}  \, \int_0 ^t
 1_{ \{ s \in [0,t] \, :\,  |B ^* e^{s
A^*} y| \ge c_0 \} }  \Big | \frac{ B^* e^{s A^*}  y }
 {|y| \, M_t } \Big|^{\alpha} ds \,
 \\ & \ge \,
  c_1 C \,  |y|^{\alpha} \, M_t ^{\alpha}  \, \int_0 ^t
   1_{ \{ s \in [0,t] \, :\,  |B ^* e^{s
A^*} y| \ge c_0 \} } \Big | \frac{ B^* e^{s A^*}  y }
 {|y| \, M_t } \Big|^{2} ds.
\eal
$$
 Let us recall  that
   the  rank condition \eqref{rank} is equivalent to the existence of
 $C_t >0$ such that, for any $u \in \R^n$,
  $\int_0 ^t  | B^* e^{s A^*}  u|^2 \, ds$
 $\ge C_t |u|^2$ (see  \cite{Za}).
Moreover
%,
% if $|y| \le 1$,
 $$
\int_0 ^t 1_{ \{ s \, :\, |B^* e^{s A^*} y| \le  c_0 \} } \, \Big
| \frac{ B^* e^{s A^*}  y }
 {|y| \, M_t } \Big|^{2} ds \le  \frac{ c_0^2 t} {|y|^2 \, M_t^2 }
 $$
%%%%%%%%%%%%%%%%
  This implies that, for any $y \in \R^n$, with $|y| \ge 1$,
 $$ \bal
& c_1 C \,  |y|^{\alpha} \, M_t ^{\alpha}  \, \int_0 ^t
   1_{ \{ s \, :\,  |B ^* e^{s
A^*} y| \ge c_0 \} } \Big | \frac{ B^* e^{s A^*}  y }
 {|y| \, M_t } \Big|^{2} ds
\\
& \ge  c_1 C C_t  \, |y|^{\alpha} \, M_t ^{\alpha -2 } - c_1 C \,
|y|^{\alpha} \, M_t ^{\alpha}  \, \int_0 ^t
   1_{ \{ s \, :\,  |B ^* e^{s
A^*} y| \le  c_0 \} } \Big | \frac{ B^* e^{s A^*}  y }
 {|y| \, M_t } \Big|^{2} ds
\\
& \ge c_1 C C_t  \, |y|^{\alpha} \, M_t ^{\alpha -2 } -
 c_1 C \,
|y|^{\alpha} \, M_t ^{\alpha} \frac{ c_0^2 t} {|y|^2 \, M_t^2 }.
\eal
$$
We get, for any $y \in \R^n, |y| \ge 1,$
$$
\int_0 ^t ds
   \int_{\R^d} \big (1 - \cos  (\lan B^* e^{s A^*}
    y, z \ran) \big)   \nu
   (dz) \ge c_1 C  C_t  M_t ^{\alpha -2 }\, |y|^{\alpha}   -
 c_1 C c_0^2 t M_t ^{\alpha -2 }.
 $$
 The assertion \eqref{df} is proved.

%%%%%%%%%%%%%%
Finally,   by the Fourier inversion formula,
\begin{equation}\label{pt}
p_t (y) = \frac{1}{(2 \pi)^n}  \int_{\R^n} e^{- i \lan y, h\ran}
\exp { \Big(- \int_0 ^t \psi (B^* e^{s A^*}  h)  ds \Big)} dh,
\;\; y \in \R^n,
\end{equation}
 is the density of $\mu_t$.
  Differentiating under the integral sign, we get easily the
 assertion. The proof is complete.\qed
 %}
 %%%%%%%%%%%%%

\begin{remark} \label{open} {\em It follows from Theorem \ref{t} that
 for any  Borel function  $f$ {\it with compact support}
  one has  $P_t f \in C^{\infty}_b (\R^n)$, for any $t>0$ (i.e.,
   $P_t f \in C^{\infty}(\R^n)$ with all bounded derivatives of
   any order)
  where
  $
P_t f (x) =  \int_{\R^n } f (  z ) p_t (z - e^{tA} x) dz.
  $
   We do not know if this regularizing effect
    holds for all   $f \in B_b (\R^n)$ as we are unable to show that
for a given  multi-index $\beta$ the
 partial derivative $D^{\beta}p_t$ is integrable on $\R^n$.
}
 \end{remark}

{{\vskip 1mm \noindent}} {\bf Acknowledgment } The authors thank
the referee
 for the careful reading  of the original manuscript, and for giving useful
 comments.

\end{document}